\documentclass{amsart}
\def \mce {{\mathcal E}}

\renewcommand{\Re}{\operatorname{\rm Re}\nolimits}
\renewcommand{\Im}{\operatorname{\rm Im}\nolimits}

\def \proj {\Pi}

\newtheorem{thm}{Theorem}[section]
\newtheorem{lemma}{Lemma}[section]

\newtheorem{cor}{Corollary}[section]

\theoremstyle{definition}

\theoremstyle{remark}

\def\squarebox#1{\hbox to #1{\hfill\vbox to #1{\vfill}}}

\def \vol {\operatorname{vol}}

\def \supp {\operatorname{supp}}
\def \card {\operatorname{card}}

\def \Vol {\operatorname{Vol}}

\def \comp {\operatorname{comp}}

\def \S1 {\mbox{S}^{1}}

\def \loc {\operatorname {loc}}

\def \Real {{\mathbb R}}
\def \Sphere {\mathbb{S}}
\def \Complex {\mathbb{C}}
\def \Natural {{\mathbb N}}

\def \mce {{\mathcal E}}

   \title [Potential scattering on cylinders]
{Asymptotics for a resonance-counting function
 for potential scattering on cylinders}
   \author { T. Christiansen}
\thanks{Partially supported by NSF grant DMS 0088922.}
\thanks{e-mail address: {\tt tjc@math.missouri.edu}}
\begin{document}
\begin{abstract} We study certain
resonance-counting functions for potential scattering
on infinite cylinders or half-cylinders.
Under certain conditions on the potential, we obtain
asymptotics of the 
counting functions, with an explicit formula for the constant appearing
in the leading term.
\end{abstract}

\maketitle
\section{Introduction}

We study potential scattering on infinite cylinders and half-cylinders.  
In particular, we give some sharp upper bounds and some asymptotics for 
resonance-counting functions in this setting.

Let $X=(-\infty,\infty)\times Y$, or $[0,\infty)\times Y$, where 
$Y$ is a smooth compact, connected manifold, with or without boundary.
We consider the product metric
$$(dx)^2+g_Y,$$
where $g_Y$ is a smooth metric on $Y$. 
Let $\Delta$ be the Laplacian on $X$, with
Dirichlet or Neumann boundary conditions
if $X$ has a boundary.  We consider operators $\Delta +V$,
where $V\in L^{\infty}_{\comp}(X;\Complex)$.

Let $\Delta_Y$ be the Laplacian on $Y$, with boundary conditions if
necessary, and
 let $\{ \sigma_j^{2}\}$, 
$\sigma_1^2\leq \sigma_2^2\leq \sigma_3^2\leq ...$ be 
the set of all eigenvalues of $\Delta_Y$, repeated according to
their multiplicity, and 
let $\nu_1^2< \nu^2_2 <\nu^2_3 <...$ be the {\em distinct}
eigenvalues of $\Delta_{Y}$.
Then the resolvent of the Laplacian $\Delta$ on $X$, or of $\Delta +V$, for
$V\in L_{\comp}^{\infty}(X)$,
 has a meromorphic continuation to the 
Riemann surface $\hat{Z}$ on which $r_j(z)=(z-\nu_j^2)^{1/2}$ is a 
single-valued function for all $j$  \cite{guce, tapsit}.  Thus
the resonances, poles of the meromorphic continuation of the 
resolvent, are points in $\hat{Z}$.  
In many settings, resonances correspond to waves which eventually decay.  
Additionally, they are in many ways analogous to eigenvalues.  Because of 
this, they have been widely studied-- see \cite{vodevsurvey, zwsurvey,
 zwqr} for an
introduction to resonances and for further references.

Here we study a simple
case of scattering on manifolds with infinite cylindrical ends.   The spectral
and scattering theory of such manifolds exhibits some characteristics one
expects both
from one-dimensional scattering and from $n$-dimensional spectral
theory (if $\dim X=n$).  The resonance-counting functions we consider here
demonstrate the one-dimensional nature of the scattering.  Evidence of 
the $n$-dimensional nature can be seen, for example, in the Weyl asymptotics
or in the maximal rate of growth of the 
eigenvalue-counting function \cite{ch-zw1, parnovski}.  It also appears in
some resonance-counting functions, e.g. \cite{cep}.

For 
$z\in \Complex \setminus [\nu^2_1, \infty)$, 
$R_V(z)=(\Delta+V-z)^{-1}$ is bounded
on $L^2(X)$ except, perhaps, for a (perhaps infinite, if 
$V$ is complex-valued) set of points corresponding to eigenvalues.
Considered as a map from $L^2_{\comp}(X)$ to $H^2_{\loc}(X)$, $R_V$ 
 has a meromorphic
continuation to the Riemann surface $\hat{Z}$ described earlier.   
Let 
$r_j(z)=(z-\nu_j^2)^{1/2}$ and let $\tilde{r}_k(z)=r_j(z)$ if 
$\sigma_k^2=\nu_j^2$.

\begin{thm}\label{thm:upperbound}  Let $X= (-\infty,\infty)\times Y$
and let $V\in L^{\infty}_{\comp}(X;\Complex)$.
Fix a sheet of $\hat{Z}$, and suppose that $\Im r_{j_0}(z)<0$ on this
sheet.  
Then, there is a constant $c_{V, \mce}\geq 0$ such that
for any $\alpha>0$, 
\begin{multline*}
\#\{ z_k: z_k \; \text{is a pole of  
}\; R_V(z) \; \text{ on this sheet},\\
|r_{j_0}(z_k)|<r,\; \Im r_{j_0}(z_k)<-\alpha\} 
= c_{V,\mce}r+o_{\alpha}(r)
\end{multline*}
The constant $c_{V,\mce}$ depends on the potential 
$V$ and the sheet (indicated by $\mce$).  Moreover, 
$$c_{V,\mce}
\leq \frac{2}{\pi}\left( \max_{(x,y),\; (x',y')\in \supp V}|x-x'| \right)
\#\{l:\Im \tilde{r}_l(z)<0 \; \text {when }\; z \; 
\text{lies on this sheet}\}.$$
\end{thm}
Here, as everywhere, we count resonances with multiplicities.  The 
error term $o_{\alpha}(r)$ depends on $V$ and on the sheet as well
as on $\alpha$, of course.

We remark that this bound on the constant $c_{V,\mce}$ is sharp, 
as can easily be seen by considering a potential that depends only
on $x$, and using the results of \cite{zw1} or \cite{froese} for potential
scattering on the line.

Although Theorem \ref{thm:upperbound} gives, in some sense,
asymptotics of a resonance-counting function, it does not give meaningful lower
bounds on the size of $c_{V,\mce}$.  In some settings we are able to actually
determine $c_{V,\mce}$, but we need some additional
conditions on $V$.

Let $\{\phi_j\}$ be an orthonormal set of eigenfunctions of 
$\Delta_Y$ associated with
$\sigma_j^2$.  
By translating if necessary, we can,
in the case of the full cylinder, arrange that for 
some $b\in \Real$, the support of $V$ is
contained in $[-b,b]\times Y$, but is not contained in the product of 
any smaller interval with $Y$.         
\begin{thm} \label{thm:asymptotics}
 Let $X=(-\infty,\infty)\times Y$ and suppose that the support
of $V$ is contained in $[-b,b]\times Y$ and the interval $[-b,b]$ cannot
be replaced by a smaller one.
Restrict ourselves
to a sheet of $\hat{Z}$ with $\Im r_j(z)<0$ if and only if $j=j_0$.
Suppose that $\nu_{j_0}^2$ is a simple eigenvalue of 
$\Delta_Y$, with $\nu_{j_0}^2=\sigma_{l_0}^2,$ and
that 
\begin{equation}\label{eq:nondegeneracy}
C|V_{l_0l_0}(x)|=C|\int_{Y}V(x,y)|\phi_{l_0}(y)|^2 d\vol_Y|
\geq |V(x,y)|,\; \text{for}\; |x-b|<\epsilon,\; |x+b|<\epsilon
\end{equation}
for some $C,\epsilon>0$.
Then,
for any $\alpha>0$, 
\begin{multline*}
\#\{ z_k: z_k\; \text{is a pole of}\; R_V(z)\; \text{ on this sheet},
|r_{j_0}(z_k)|<r,\; \Im r_{j_0}(z_k)<-\alpha\} \\
= \frac{4}{\pi}b 
r +o_{\alpha}(r).
\end{multline*}
\end{thm}
In Section \ref{s:pta} we give an example of a 
nontrivial complex-valued potential for
which (\ref{eq:nondegeneracy}) is not satisfied and for which the 
conclusion of the theorem does not hold.  Moreover, for this potential
$c_{V,\mce}=0$ for at least one (non-physical) sheet.  
This gives an example of some behaviour which is even asymptotically
 truly different 
from that demonstrated by scattering by the family of potentials
$V(x)$.  Moreover, this means that
potential scattering on cylinders provides an example of a setting in 
which even the order of growth of a resonance-counting function may vary 
depending on the potential.

In Section \ref{s:pta} we prove a theorem which gives another situation
in which we can determine $c_{V,\mce}$.
In Section \ref{s:halfcylinder} we give some analogous results for
potential scattering on half-cylinders.

Scattering on cylinders bears some resemblance to potential scattering on
the line.  On the line, the distribution of resonances has been studied
in \cite{froese,regge,zw1}.   The 
complicated nature of $\hat{Z}$ makes more
difficult the question
of bounding the number of resonances in the cylindrical end setting.
 Earlier
results on 
resonances for manifolds with cylindrical ends include \cite{a-p-v, cep,
d-e-m, edward}, and references.
For general scattering theory on manifolds
with cylindrical ends, references include \cite{guce, tapsit}.

\section{Preliminaries and notation}

Let $r_j(z)=(z-\nu_j^2)^{1/2}$ and 
identify the physical sheet of 
$\hat{Z}$ as being the part of $\hat{Z}$ on which $\Im r_j(z)>0$ for
all $j$ and all 
$z$ and on which $R_V(z)$ is bounded
on $L^2(X)$ for all but a discrete set of $z$.
  Other sheets will be identified,
when necessary, by indicating for which values of $j$ $\Im r_j(z)<0$. 
Each sheet can be identified with  $\Complex \setminus [\nu_1^2,\infty)$.
With this language, there are points in $\hat{Z}$ which belong to no sheet
but which belong to the boundary of the closure of two sheets, and the 
ramification points, which correspond to $\{\nu_j^2\}$ and belong to 
the closure of four sheets (except for ramification points
corresponding to $\nu_1^2$).  We note that sheets that meet the physical sheet
are characterized by the existence of a $J\in \Natural $ such that
$$\Im r_j(z)<0 \; \text{for all } z \;\text{on that sheet if and only if} \;
j \leq J.$$

 We can associate to a fixed sheet of $\hat{Z}$ a set 
$\mce \subset \Natural$,
$$\mce = \{ j: \Im r_j(z)<0\; \text{on this sheet}\}.$$  We shall call 
$\mce $ the {\em labeling set}.
Let
$$\tilde{\mce}=\{ l\in \Natural : \sigma_l^2=\nu_j^2\;
 \text{for some } \; j\in \mce\}.
$$ 
Let $\{\phi_j\}$ be an orthonormal set of eigenfunctions of $\Delta_Y$
associated with $\{\sigma_j^2\}$.

In general, we shall use $z$ to stand for a point in $\hat{Z}$ and 
$\proj(z)$ to represent its projection to $\Complex$.  For 
$w\in \Real^m$, $\langle w \rangle =(1+|w|^2)^{1/2}$.  We will denote
by $C$ a constant whose value may change from line to line.

Next we recall some results and language of complex analysis,
e.g. \cite{levin},
and recall a theorem we shall need on the distribution of zeros of functions
which are ``good'' in a half-plane.

We shall often
work with functions that are holomorphic not in the whole
plane but are holomorphic within an angle $(\theta_1, \theta_2)$.
A function $F$ holomorphic in an angle $(\theta_1,\theta_2)$ is of order
$\rho$ there
 if $$\overline{\lim}_{r\rightarrow \infty}\frac{\ln \ln (\sup_{\theta 
\in \theta_1,\theta_2}|F(re^{i\theta})|)}{r}=\rho.$$
A function of order $\rho$ in the angle $(\theta_1,\theta_2)$ is of type
$\tau$ there if 
$$\overline{\lim} _{r\rightarrow \infty} \frac{\ln 
\sup_{\theta\in (\theta_1,\theta_2)}|F(re^{i\theta})|}{r^{\rho}}=\tau.$$
 A function of order
$1$ and type $\tau<\infty$ (in an angle
$(\theta_1,\theta_2)$) is said to be of exponential type there.  Of course,
$\rho$ and $\tau$ can depend on $\theta_1$ and $ \theta_2$.

The indicator of a function $F$ holomorphic in an angle $\theta_1 \leq \arg 
\zeta
\leq \theta_2$
and of order $\rho$ is 
$$h_F(\theta)=\overline{\lim}_{r\rightarrow \infty}\frac{\ln |F(re^{i\theta})|}
{r^{\rho}}.$$
A function $F$ is of completely regular growth within the angle
$(\theta_1,\theta_2)$ if 
$$\lim_{r\rightarrow \infty \atop r\not \in E}\frac{\ln |F(re^{i\theta})|}
{r^{\rho}}=h_F(\theta)
$$
where the set $E\subset \Real_+$ is of zero relative measure and the 
convergence is uniform for $\theta\in (\theta_1,\theta_2)$.

We shall abuse notation slightly and also use the language above for a function
that is holomorphic for $\theta_1 \leq \arg \zeta \leq \theta_2$ and 
$\zeta$ outside of a compact set.

For a function $f$ defined in the lower half plane, let 
$n_{f-}(r)$ be the number of zeros of $f$, counted with multiplicity, that
lie in the lower half-plane and have norm less than $r$.
\begin{thm}\label{thm:complexanal}
Suppose $f(\zeta)$
 is holomorphic in the closed lower half plane $\Im \zeta \leq 0$,
\begin{equation*}\label{as1}
|f(\zeta )|\leq Ce^{C|\zeta |}
\end{equation*}
 there, $f(0)=1$,
\begin{equation*}\label{as3}
\left |\int_{-\infty}^{\infty}\frac{d[\arg f(t)]}{dt}dt\right|<\infty
\end{equation*}
and 
\begin{equation*}\label{as4}
\left|\int_{-\infty}^{\infty}\frac{\ln |f(t)|}{1+t^2}dt\right|<\infty.
\end{equation*}
Then 
$$\lim_{r\rightarrow \infty}\frac{n_{f-}(r)}{r}
=\frac{1}{2\pi} \int_{\pi}^{2\pi}h_f(\varphi)d\varphi.$$
\end{thm}
The proof of this theorem can be found in \cite{steplike}.  It is an 
adaptation of arguments of \cite[Chapter III, Section 2]{levin}
and
\cite[Theorem 3, Chapter III, Section 3]{levin}.

We note, moreover, that the assumptions of Theorem \ref{thm:complexanal} 
mean that $f$ is a function of completely regular growth in the lower 
half-plane and that $h_f(\theta)=c_f|\sin \theta|$ for $\pi < \theta< 2\pi$.

\section{Proof of Theorem \ref{thm:upperbound}}\label{s:pfthmub}

As in \cite{froese}, here we find a matrix $B$ so that the poles of 
the resolvent in the region in question are included in the zeros
of $\det(I+B)$.  We study the properties of 
the matrix $B$, and then apply Theorem
\ref{thm:complexanal}.  Recall that here $X=(-\infty,\infty)\times Y$.

Let \begin{equation}
R_0(z)=
(\Delta -z)^{-1}=
\sum_{j=1}^{\infty}\frac{i}{2r_j(z)}e^{i|x-x'|r_j(z)}
\sum_{\sigma_l^2=\nu_j^2}\phi_l(y)\overline{\phi}_l(y').
\end{equation}
Then 
$$(\Delta+V-z)R_0(z)=I+VR_0(z).$$
Since $R_0(z)$ has no null space, away from the ramification points of
$\hat{Z}$, $R_V(z)$ has a pole if and only if $I+VR_0(z)$ has nontrivial
null space (and the multiplicities agree).

If $\mce \subset \Natural$ is a finite set, define 
$w_{\mce}:\hat{Z}\rightarrow \hat{Z}$ as follows.
To $z$ we may associate the set of square roots $\{r_j(z)\}$.  Then 
$w_{\mce}(z)$ may
be determined by saying it is the element of $\hat{Z}$  
associated to the set $\{ r_j(w_{\mce}(z))\}$, with
$$r_j(w_{\mce}(z))=\left\{ \begin{array}{ll} -r_j(z),\; &
 \text {if}\; j \in \mce \\
r_j(z),& \text {if}\; j\not \in \mce.
\end{array} \right.
$$  

Suppose we now restrict ourselves to consider only $z$
lying on the sheet with 
$$\Im r_j(z)<0 \; \text{if and only if}\; j\in \mce.$$
Then $w_{\mce}(z)$ lies in the physical sheet.  
Moreover,
\begin{align}
I+VR_0(z)& = \left(I+VR_0(w_{\mce}(z))\right)\left[
I+\left[I+ VR_0(w_{\mce}(z))\right]^{-1}V\left[R_0(z)-R_0(w_{\mce}(z))\right]
\right] \nonumber \\
& = \left(I+VR_0(w_{\mce}(z))\right)
\left[I+ \left[I+VR_0(w_{\mce}(z))\right]^{-1}
A_1(z)\right]
\end{align}
where $A_1(z)$ has Schwartz kernel
$$V(x,y)
\sum_{l\in \tilde{\mce}}\frac{i}{2\tilde{r}_l(z)}(e^{i\tilde{r}_l(z)(x-x')}
+ e^{-i\tilde{r}_l(z)(x-x')})\phi_l(y) \overline{\phi_l}(y').$$
If $|\Im \proj(w_{\mce}(z))|> \|\Im V(x,y)\|_{L^{\infty}}$,
then $I+VR_0(w_{\mce}(z))$ is invertible.  If we restrict ourselves to
such $z$, then, the poles of the resolvent of $\Delta +V$ are given by the
zeros of 
$$\det(I+A_2(z)),$$
where $A_2(z)$ is 
$$A_2(z)=\sum_{l\in \tilde{\mce}}\frac{i}{2\tilde{r}_l(z)}
( \varphi_{l,+}\otimes \Psi_{l,-}
+ \varphi_{l,-}\otimes \Psi_{l,+}),$$
with  
\begin{align*}
\Phi_{l\pm}(x,y,z)& = e^{\pm i\tilde{r}_l(z)x}\phi_{l}(y)\\
\varphi_{l\pm}(x,y,z)&=\left(\left(I+VR_0(w_{\mce}(z))\right)^{-1}
(V \Phi_{l\pm})
(\bullet,z)\right)(x,y)\\
\Psi_{l,\pm}(x,y,z)&=e^{\pm i \tilde{r}_l(z)x}\overline{\phi_l}(y).
\end{align*}
Here we use the notation 
$$(f\otimes g)h(x,y)=f(x,y)\int_X g(x',y') h(x',y')d\vol_{X}.
$$

One can then see that the zeros of $I+A_2(z)$ are the same as the zeros
of $I+A_2(z)\chi$, where $\chi\in L^{\infty}_{\comp}
(X)$ is one on the support of
$V$.  The zeros of $I+A_2(z)\chi$ are the same as the zeros of 
$\det(I+B(z))$, where 
\begin{equation}\label{eq:B}
B(z)=\left( \begin{array}{cc}
B_{+-}(z) & B_{--}(z)\\
B_{++}(z) & B_{-+}(z)
\end{array}\right),
\end{equation}
$B_{+\pm}=(b_{+\pm lj})_{l,j\in \tilde{\mce}}$, 
$B_{-\pm}=(b_{-\pm lj})_{l,j\in \tilde{\mce}}$,
and
\begin{align} \label{eq:b+-}
\nonumber
b_{+\mp lj}(z)& = 
\frac{i}{2\tilde{r}_l(z)}\int_X \varphi_{j+}(x,y,z)\chi(x,y)\Psi_{l\mp}
(x,y,z)d\vol_X,
\\
b_{-\mp lj}(z)& = 
\frac{i}{2\tilde{r}_l(z)}\int_X \varphi_{j-}(x,y,z)\chi(x,y)\Psi_{l\mp}
(x,y,z)d\vol_X.
\end{align}

We shall first obtain upper bounds on the entries in the matrix $B$,
and thus on $\det(I+B(z))$.  To do so, we will use the following
lemma.
\begin{lemma}\label{l:r0bound}
Let $f_{\pm}(x,z)=e^{\pm i \tilde{r}_j(z)x}$, and let $\chi_1,\; \chi_2 \in
C_c^{\infty}(X)$.  If $z$ lies on the physical sheet of $\hat{Z}$ and
$\Im \tilde{r}_j(z)=t_0>0$, then 
$$\left \| \chi_1 \frac{1}{f_{\pm}}R_0(z)f_{\pm}\chi_2
\right \|_{L^2(X)\rightarrow 
L^2(X)}
\leq \frac{C}{|\Re \tilde{r}_j(z)|^{7/12}}$$
when $|\tilde{r}_j(z)|$ is sufficiently large.  Moreover, for $\Im \tilde{r}_j
(z) \geq t_0>0$,
$$\left \| \chi_1 \frac{1}{f_{\pm}}R_0(z)f_{\pm}\chi_2
\right \|_{L^2(X)\rightarrow 
L^2(X)}
\leq \frac{C}{| \tilde{r}_j(z)|^{5/12}}$$
when $|\tilde{r}_j(z)| $ is sufficiently large.
\end{lemma}
\begin{proof}
Without loss of generality we can assume $\chi_1$ and
$\chi_2$ are independent of $y$ and 
thus it is suffices to consider, for $l\in \Natural$, 
$$\left \| \chi_1 \frac{1}{f_{\pm}}R_{0l}(z)f_{\pm}\chi_2
\right\|_{L^2(X)\rightarrow 
L^2(X)}$$
where $R_{0l}(z)$ has Schwartz kernel 
$$\frac{i}{2\tilde{r}_l(z)}e^{i\tilde{r}_l(z)|x-x'|}\phi_l(y)\overline{\phi}_l
(y').$$

The Schwartz kernel of $(f_{\pm})^{-1}R_{0l}(z)f_{\pm}$ is 
$$K_{l\pm}(x,y,x',y',z)=
\left\{ \begin{array}{ll}
\frac{i}{2\tilde{r}_l(z)}e^{i(\tilde{r}_l(z)\mp\tilde{r}_j(z))(x-x')}
\phi_l(y)\overline{\phi}_l(y') \chi_1(x)\chi_2(x'),\; & \text{if}\; 
x>x'\\\frac{i}{
2\tilde{r}_l(z)}e^{i(-\tilde{r}_l(z)\mp\tilde{r}_j(z))(x-x')}
\phi_l(y)\overline{\phi}_l(y') \chi_1(x)\chi_2(x'),\; &\text{if}\; 
x<x'.
\end{array}
\right.
$$
We shall show that when $\Im \tilde{r}_j(z)=t_0$
$$\int_X\int_X |K_{l\pm}(x,y,x',y',z)|^2 d\vol_X d\vol_X \leq
\frac{C}{|\Re \tilde{r}_j(z)|^{7/6}},$$
with constant $C$ independent of $l$, which will prove the first
part of the lemma.

First, notice that on the support of 
$\chi_1(x)\chi_2(x')$,
the exponential function in $K_{l\pm}$ is bounded 
independent of $l$.  
This is because $\Im \tilde{r}_l(z)>0$ and $|\Im \tilde{r}_j(z)(x-x')|$
is bounded  for $x\in \supp \chi_1$, $x'\in \supp \chi_2.$  Thus,
\begin{equation}\label{eq:bbd}
\|K_{l\pm}(z)\|^2_{L^2}\leq \frac{C}{|\tilde{r}_l(z)|^2}.
\end{equation}

When $\tilde{r}_l\not = \tilde{r}_j$, we may integrate by parts to see
that 
$$\| K_{l\pm}(z)\|^2_{L^2}\leq \frac{C}{|\Im( \tilde{r}_j(z)-\tilde{r_l}(z))|}
\frac{1}{|\tilde{r}_l(z)|^2}$$
so that
$$\| K_{l\pm}(z)\|^2_{L^2}\leq 
\frac{C}{|\tilde{r}_l(z)|^2} \min(1,
(|\Im( \tilde{r}_j(z)-\tilde{r_l}(z))|)^{-1}).$$

Let $\tilde{r}_j=s+it_0$.  Then if $\tilde{r}_l(z)=u+iv$, a computation 
shows that, with $g=\sigma_j^2+s^2-t_0^2-\sigma_l^2$, 
$u^2=\frac{1}{2}(g+\sqrt{g^2+4s^2t_0^2})$, and 
$v^2=\frac{1}{2}(-g+\sqrt{g^2+4s^2t_0^2})$.
If $g\leq (|s|t_0)^{7/6},$ then 
\begin{align}\label{eq:vbound}
v^2 & \geq \frac{1}{2}\left( -(|s|t_0)^{7/6}+\sqrt{(|s|t_0)^{7/3}+
4(|s|t_0)^{2}}\right)\nonumber \\
& = (|s|t_0)^{5/6}+O((|s|t_0)^{1/2}).
\end{align}
Then $$\|K_{l\pm}(z)\|^2_{L^2}\leq \frac{C}{|\tilde{r}_l(z)|^2|v-t_0|}
\leq \frac{C}{(|s|t_0)^{5/6}(|s|t_0)^{5/12}}\leq \frac{C}{|s|^{5/4}}$$
when $|s|$ is sufficiently large and $\Im \tilde{r}_j(z)=t_0$.

If, on the other hand, $g\geq (|s|t_0)^{7/6}$, then we use 
\begin{equation}\label{eq:ubound}
u^2=\frac{1}{2}(g+(\sqrt{g^2+4(st_0)^2})\geq g \geq (|s|t_0)^{7/6}
\end{equation}
and 
$$\|K_{l\pm}(z)\|_{L^2}^2\leq \frac{C}{|\tilde{r}_l(z)|^2}
\leq \frac{C}{u^2}\leq \frac{C}{(|s|t_0)^{7/6}}.$$  This finishes the proof
of the first part of the lemma.

To prove the second part of the lemma, first notice that
if $\tilde{r}_j(z)=s+it$ and $|s|<1$, then 
$\frac{1}{|\tilde{r}_l(z)|^2}\leq \frac{C}{t^2}$ when $t$ is sufficiently
large, so that $$\|K_{l\pm}(z)\|_{L^2}^2\leq \frac{C}{t^2}$$ in this region.
On the other hand, if $|s|\geq 1$, the inequalities (\ref{eq:bbd}),
(\ref{eq:vbound}) and 
(\ref{eq:ubound}) together show that when $t\geq t_0$, 
$$\|K_{l\pm}(z)\|^2_{L^2} \leq \frac{C}{|s+it|^{5/6}}.$$

\end{proof}

Fix $j_0\in \mce$.  We shall eventually
use $k=r_{j_0}(z)$ to identify our fixed
sheet of $\hat{Z}$ (corresponding to $\mce$) with the lower half plane.
However, we shall continue to use $z$ as a coordinate as well, when
it is more convenient.  In any case, we restrict ourselves to one fixed
sheet.

\begin{lemma}\label{l:mixedbd}  Fix a sheet of $\hat{Z}$ with corresponding
labeling set $\mce$ and let $j_0\in \mce$, $l,j\in \tilde{\mce}$.
If $-\Im r_{j_0}(z)\geq\alpha>0$, then for $|r_{j_0}(z)|$ sufficiently
large (depending on $\alpha$),  
$$|b_{+-lj}(z)|\leq \frac{C}{|\tilde{r}_l(z)|},\; |b_{-+lj}(z)|\leq \frac
{C}{|\tilde{r}_l(z)|}.$$
\end{lemma}
\begin{proof}
First we show that 
in this region, for $j\in \tilde{\mce}$ and $\chi \in L^{\infty}_{\comp}(X)$,
\begin{equation}\label{eq:conjresbd}
\| e^{\mp i\tilde{r}_j(z)x}(I+VR_0(w_{\mce}(z)))^{-1}
\chi e^{\pm i\tilde{r}_j(z)x}
\| \leq C
\end{equation}
when $|r_{j_0}(z)|$ is sufficiently large.

When $|r_{j_0}(z)|$ is sufficiently large, and $\tilde{\chi}\in 
L^{\infty}_{\comp}(X)$ is one on the support of $V$, 
\begin{align}\label{eq:mixedsignbd}& 
\| \tilde{\chi}e^{\mp i\tilde{r}_j(z)x}(I+VR_0(w_{\mce}(z)))^{-1}\chi
e^{\pm i \tilde{r}_j(z)x}\| \\ \nonumber &
= \| \sum_{m=0}^{\infty} e^{\mp i\tilde{r}_j(z)x} (-1)^m
(VR_0(w_{\mce}(z))\tilde{\chi})^m \chi  e^{\pm i \tilde{r}_j(z)x} \|\\  \nonumber &
= \| \sum_{m=0}^{\infty} (-1)^m
(e^{\mp i\tilde{r}_j(z)x} 
VR_0(w_{\mce}(z))\tilde{\chi}e^{\pm i\tilde{r}_j(z)x})^m \chi  \|\ \\ \nonumber &
\leq C 
\end{align}
where we are using Lemma \ref{l:r0bound}.
Using this estimate and the definition 
of $b_{+-lj}$, $b_{-+lj}$, we obtain the 
desired estimates.
\end{proof}

We shall need the following bound on the $b_{++lj}(z)$ and $b_{--lj}(z)$.
\begin{lemma}\label{l:samesignbd}  Fix a sheet of $\hat{Z}$ with corresponding
labeling set $\mce$, and let $j_0\in \mce$.
If $\Im r_{j_0}(z)\leq -\alpha<0$, $l,j\in \tilde{\mce},$
and $\supp (V) \subset [-\beta,\gamma]$,  then for $|r_{j_0}(z)|$ sufficiently
large (depending on $\alpha$),  
$$|b_{++lj}(z)|\leq \frac{Ce^{2\gamma|\Im \tilde{r}_j(z)|}}
{|\tilde{r}_l(z)|},\; 
|b_{--lj}(z)|\leq \frac
{Ce^{2\beta|\Im \tilde{r}_j(z)|}}{|\tilde{r}_l(z)|}.$$
\end{lemma}
\begin{proof}
We give the proof for $b_{++lj}$.    Note that if $\supp f \subset \supp V$,
then $\supp(I+VR_0(w_{\mce}(z)))^{-1}f\subset \supp V$.
Recall that 
$\Phi_{j\pm}(x,y,z) = e^{\pm i\tilde{r}_j(z)x}\phi_{j}(y).$
Then as in (\ref{eq:mixedsignbd}), we obtain
that 
$$
\| \varphi_{j+}\|=
\|(I+ VR_0(w_{\mce}(z)))^{-1}V\Phi_{j+}\| \leq C e^{\gamma
|\Im \tilde{r}_j(z)|}.
$$
Using this bound and the remark about the support properties of 
$(I+VR_0(\omega_{\mce}(z)))^{-1}$, 
\begin{align*}
|b_{++lj}(z)| & 
=\left|\frac{1}{2\tilde{r}_l(z)}\int_{X} \varphi_{j+}(x,y,z)
\chi(x,y)\Psi_{l+}(x,y,z)
d\vol_X\right| \\ & 
\leq \frac{C}{|\tilde{r}_l(z)|}e^{2\gamma |\Im \tilde{r}_j (z)|}
\end{align*}
for $|r_{j_0}(z)|$ sufficiently large, $\Im r_{j_0}(z)\leq -\alpha$.
For the last inequality we have also used that $\tilde{r}_i(z)
\rightarrow\tilde{r}_j(z)$ as $\proj(z)\rightarrow \infty$.

A similar argument yields the proof of the bound for $b_{--lj}(z)$.
\end{proof}

\begin{proof}[Proof of Theorem \ref{thm:upperbound}]
We use the coordinate $k=r_{j_0}(z)$ to identify our fixed sheet with
the lower half plane.  Let $g_1(k)=\det(I+B(z(k)))$, where $\proj(z(k))=
k^2+\nu_{j_0}^2$ and $z$ lies on our sheet.
Here $B(z)$ is as defined in (\ref{eq:B}) and (\ref{eq:b+-}).  
Then $g_1(k)$ has
at most a finite number of poles, $k_{1}, k_{2},...,k_{m_{\alpha}}$,
listed with multiplicity, in $
\Im k
\leq -\alpha$.
Let 
$$g_2(k)=g_1(k)(k-k_{1})(k-k_{2})\cdot\cdot\cdot
(k-k_{m_{\alpha}})$$
and, if $g_2(-i\alpha)\not = 0$, let
$$g_3(k)=\frac{g_2(k)}{g_2(-i\alpha)}.$$
If $g_2(-i\alpha)=0$, let 
$$g_3(k)=\frac{g_2(k)l!}{(k+i\alpha)^lg_2^{(l)}(-i\alpha)}$$
where $l$ is chosen so that $g_2^{(m)}(-i\alpha)=0$ if $m<l$ but $g_2^{(l)}
(-i\alpha)\not = 0$.
Then Lemmas \ref{l:mixedbd} and \ref{l:samesignbd} show that
the hypotheses of Theorem \ref{thm:complexanal} are satisfied for
$g_4(k)=g_3(k-i\alpha)$, with 
$$|h_{g_4}(\varphi)|\leq 2\left(\sup_{(x,y),\;(x',y')
\in \supp V}|x-x'|\right)\card(\tilde{\mce})
|\sin \varphi|
.$$
Recalling that, except possibly for a finite number, the zeros of
$g_3(k)$ correspond to the poles of $R_V(z)$ in this region, an application
of Theorem \ref{thm:complexanal} finishes the proof.
\end{proof}

\section{Determining $c_{V,\mce}$ and a counterexample}
\label{s:pta}

In this section we prove Theorem \ref{thm:asymptotics}, give a 
counterexample, and give another example of a setting in which 
$c_{V,\mce}$ can be determined.

We shall need the following lemma, which is Lemma 4.1 of \cite{froese}.
\begin{lemma}\label{l:froese}
Suppose $v\in L^{\infty}(\Real)$ has compact support contained in $[-1,1]$,
but in no smaller interval.   Suppose $f(x,k)$ is analytic for $k$ in
the lower half plane, and for real $k$ we have $f(x,k)\in L^2([-1,1]dx,\Real
dk)$.  Then $\int e^{\pm ikx}v(x)(1-f(x,k))dx$ has exponential type at least 
$1$ for $k$ in the lower half plane.
\end{lemma}

In the next lemma, we use $k=\tilde{r}_l(z)$ as a coordinate, and, fixing a 
sheet of $\hat{Z}$, let $z(k)$ be the corresponding point on $\hat{Z}$.
\begin{lemma}\label{l:mintype}
 Let $X=(-\infty,\infty)\times Y$ and suppose that the support
of $V$ is contained in $[-b,b]\times Y$ and the interval $[-b,b]$ cannot
be replaced by a smaller one.
Suppose
that 
\begin{equation*}
C|V_{ll}(x)|=C|\int_{Y}V(x,y)|\phi_{l}(y)|^2 d\vol_Y|
\geq |V(x,y)|,\; \text{for}\; |x-b|<\epsilon,\; |x+b|<\epsilon
\end{equation*}
for some $C,\epsilon>0$.  Fix a sheet of $\hat{Z}$ on which
$\Im \tilde{r}_l(z)<0$, and choose $\alpha$ so that there
are no poles of $b_{++ll},\; b_{--ll}$
on this sheet with $\Im r_l(z)\leq -\alpha$.  
Then $b_{++ll}(z(k)),\; b_{--ll}(z(k))$
are functions of type at least $2b$ for 
the half-plane $\Im k \leq -\alpha$, $k=\tilde{r}_l (z)$.
\end{lemma}
\begin{proof}
We give the proof for $b_{++ll}$, as the proof for $b_{--ll}$ is 
similar.

Let $g(k,x)=e^{ikx}$ and $\mce$ be the labeling set
associated to our fixed sheet of $\hat{Z}$.  Let 
\begin{align*}
f_1(x,y,k)& =\overline{\phi}_l(y)\frac{1}{V(x,y)}
\left[\frac{1}{g}\left[ I-[I+VR_0(\omega_{\mce}(z(k)))]^{-1}\right]
V \Phi_{l+}(\cdot, z(k))\right](x,y)
\\ & =\overline{\phi}_l(y) 
\sum_{m=1}^{\infty}(-1)^m \left[\left((g)^{-1}R_0(\omega_{\mce}(z(k)))
Vg\right)^m \phi_l\right](x,y)
\end{align*}
where the second equality holds when $| k|$ is sufficiently large.
Then 
\begin{align*}
b_{++ll}(z(k))& =\frac{i}{2 k}\int 
e^{2ikx}V(x,y)(|\phi_l|^2(y)-f_1(x,y,k))d\vol_X.
\end{align*}

Let 
\begin{equation*}
\chi_{\epsilon}(x)=\left\{
\begin{array}{ll}
0, \;&  \text{if}\; |x|<b-\epsilon\; \text{or}\; |x|>b\\
1, & \text{if}\; b-\epsilon\leq |x|\leq b.
\end{array}
\right.
\end{equation*}
Let 
$$v(x)=\int_Y V(x,y)|\phi_l|^2(y) d\vol_Y=V_{ll}(x),$$ and 
$$f(x,k)=
\frac{1}{V_{ll}(x)}\chi_{\epsilon}(x)
\int_Y V(x,y)f_1(k,x,y)d\vol_Y.$$
Note that
$$b_{++ll}(z(k)) =\frac{i}{2 k}\int 
e^{2ikx}v(x)(1-f(x,k))dx 
-\int_X e^{2ikx}(1-\chi_{\epsilon})V(x,y)f_1(x,y,k)d\vol_X.$$
Using (\ref{eq:conjresbd}) and the support properties of 
$V(1-\chi_{\epsilon})$, 
the last term on the right is of type at most $2b-2\epsilon$, and so we need only show that
the first integral on the right is of type at least $2b$.  To do this, we
will apply Lemma \ref{l:froese} to $b_{++ll}(z(k+i\alpha)).$

We must show that $f(x,k)\in L^2([-b,b]dx,\Real dk)$ when $\Im k =
-\alpha$.
We have
\begin{align*}
\int |f(x,k)|^2dx &= \int_{|x|\leq b} |V_{ll}(x)|^{-2}\chi_{\epsilon}
(x)\left|\int_Y
V(x,y)f_1(k,x,y)d\vol_Y\right|^2dx\\ &
\leq   \int_{|x|\leq b}\int_{Y}|V_{ll}(x)|^{-2}\chi_{\epsilon}(x)
|V(x,y)|^2d\vol_Y \int_Y |f_1(k,x,y)|^2 d\vol_Y dx\\
&  \leq C \int_X |f_1(k,x,y)|^2d\vol_X.
\end{align*}
By Lemma \ref{l:r0bound}, when $|\Re k|$ is sufficiently large, this is 
bounded by $C|\Re k|^{-7/6}.$  When $|\Re k|$ is in a compact set (with
$\Im k =-\alpha$), it is enough to note that 
$\int|f_1(k,x,y)|^2d\Vol_X$ is bounded, so that 
$f(x,k)\in L^2([-b,b]dx,\Real dk).$
Then, applying Lemma \ref{l:froese} after appropriately
rescaling, we finish the proof.
\end{proof}

\begin{proof}
[Proof of Theorem \ref{thm:asymptotics}]  We use 
$k=r_{j_0}(z)=\tilde{r}_{l_0}(z)$ as the coordinate.
The simplicity of $\nu_{j_0}^2$ as an eigenvalue of $\Delta_{Y   } $
means that the matrix $B$ is a $2\times 2$ matrix
$$B=\left( \begin{array}{cc}
b_{+-l_0l_0}& b_{--l_0l_0}\\
b_{++l_0l_0}& b_{-+l_0l_0}
\end{array}
\right).$$
Thus $\det(I+B)(z(k))=
[(1+b_{+-l_0l_0})(1+b_{-+l_0l_0})-b_{--l_0l_0}b_{++l_0l_0}](z(k))
=\varphi_1(k).$

Suppose first that
 $\varphi_1(k)$ has no poles in the region $\Im k\leq -\alpha$.  
If $\varphi_1(-i\alpha)\not = 0$, let
$$\varphi_2(k)=\frac{\varphi_1(k)}{\varphi_1(-i\alpha)}.$$
If $\varphi_1(-i\alpha)=0$, let 
$$\varphi_2(k)=\frac{\varphi_1(k)l!}{(k+i\alpha)^l\varphi_1^{(l)}(-i\alpha)}$$
where $l$ is chosen so that $\varphi_1^{(m)}(-i\alpha)=0$ if $m<l$ but 
$\varphi_1^{(l)}
(-i\alpha)\not = 0$.  

Note that by Lemmas \ref{l:mixedbd}
and \ref{l:samesignbd}, 
for $s \in \Real$, $\varphi_2(s-i\alpha)=c_0(1+O(|s|^{-1}))$ when $|s|
\rightarrow \infty$, for some nonzero constant $c_0$.  Moreover, by 
Lemmas \ref{l:mixedbd}, \ref{l:samesignbd}, and \ref{l:mintype},
$\varphi_2(k)$ is a function of type $4b$ in the half plane
$\Im k \leq -\alpha $.  Then applying
Theorem \ref{thm:complexanal} to $\varphi_2(k)$ in 
the half-plane $\Im k \leq -\alpha$, we obtain the result.

If $\varphi_1(k) $ has poles in the region $\Im k\leq - \alpha$, they 
can be handled in the same manner as in the proof of Theorem 
\ref{thm:upperbound}.
\end{proof}
 
We give a counterexample for Theorem \ref{thm:asymptotics}.
Let $X=\Real \times \Sphere^1$, and let $V(x,y)=V_1(x)e^{imy}$ with
$V_1(x)\in L^{\infty}_{\comp}(\Real)$ nontrivial and $m>0$ an integer.
Then on the sheet with $\Im r_j(z)<0$ if and only if $j=0$ there are
no resonances.  To see this,
restrict $z$ to this sheet of $\hat{Z}$ and 
note that since for integers $n$
$R_0$ commutes with projection onto the span of
$e^{iny}$, 
$$(VR_0(\omega_{\{0\}}(x)))^jV\Phi_{0,\pm}= u_{\pm,j}(x)e^{im(j+1)y}$$
for $j=0,1,2,...$ .
Moreover, 
$$(I+VR_0(\omega_{\{0\}}(z)))^{-1}V= 
\sum_{j=0}^{\infty}(-1)^j(VR_0(\omega_{\{0\}}(z)))^jV$$
 when $-\Im r_0(z)
\geq \alpha >0$ and $|r_0(z)|$ is sufficiently large.  Therefore, using 
(\ref{eq:b+-}), we see that $b_{\pm +00}(z)=0=b_{\pm -00}(z)$ when
$-\Im r_0(z)
\geq \alpha >0$ and $|r_0(z)|$ is sufficiently large.  By analytic
continuation, $b_{\pm +00}(z)=0=b_{\pm -00}(z)$ for all $z\in \hat{Z}$.  Thus,
by the discussion of Section \ref{s:pfthmub}, there are no resonances on
the sheet with $\Im r_j(z)<0$ if and only if $j=0$. 

In Theorem \ref{thm:asymptotics} we used some knowledge of the potential 
near the boundary of its support to allow us to find $c_{V,\mce}$.  In 
the following theorem we again make use of the fact that the potential 
is ``controlled'' near the boundary of its support.
\begin{thm}
Suppose for some potential $V_0\in L^{\infty}_{\comp}(X;\Complex)$, with 
$\supp V_0\subset [-b_0,b_0]\times Y$ and for some sheet of $\hat{Z}$ 
with corresponding
labeling set $\mce\ni j_0$, we have
\begin{multline*}
\#\{ z_k: z_k \; \text{is a pole of}\; R_{V_0}(z)\;  
 \text{on this sheet},\\
|r_{j_0}(z_k)|<r,\; \Im r_{j_0}(z_k)<-\alpha\} 
= \frac{4b_0}{\pi}\#\{l: \; l\in \tilde{\mce}\}  r+o_{\alpha}(r)
\end{multline*}
for some $\alpha>0$.
Suppose in addition $W\in L^{\infty}_{\comp}(X;\Complex)$ with
$\supp W\subset [-b_0+\epsilon, b_0-\epsilon]\times Y$ for some $\epsilon >0$.
Then 
\begin{multline*}
\#\{ z_k: z_k \; \text{is a pole of
}\; R_{V_0+W}(z) \; \text{on this sheet},\\
|r_{j_0}(z_k)|<r,\; \Im r_{j_0}(z_k)<-\alpha\} 
= \frac{4b_0}{\pi}\#\{l: \; l\in \tilde{\mce}\}  r+o_{\alpha}(r).
\end{multline*}
\end{thm}
That is, if the resonance-counting function for $\Delta+V_0$ has maximal
growth rate, so does that for $\Delta+V_0+W$.
\begin{proof}
In the proof of this theorem, we will add a superscript to the matrix 
$B$ from 
Section \ref{s:pfthmub} and its entries to indicate to which
potential it is associated.  That is, when $|r_{j_0}(z)|$ is sufficiently
large, the poles of the resolvent of
 $\Delta +V$ correspond to the zeros of $\det(I+B^V(z))$
and likewise for $V_0$.

In this proof, as previously, we shall sometimes use as coordinate on 
our sheet $k=r_{j_0}(z)$, and then $z(k)$ is the corresponding point 
on our sheet.

Let $V=V_0+W$.
We shall show that $B^V(z)=B^{V_0}(z)+D(z)$, with the entries $d_{lj}(z)$
of $D(z)$ satisfying
\begin{equation}\label{eq:dbounds}
|d_{lj}(z)|\leq \frac{C}{|\tilde{r}_{l}(z)|}
e^{(2b_0-\epsilon)|\Im \tilde{r}_j(z)|}.
\end{equation}
Because of the assumption on the distribution of resonances for 
$\Delta +V_0$, $\det(I+B^{V_0}(z(k)))$ is of type 
$4b_0 \#\{l: \; l\in \tilde{\mce}\}$ in
$\Im k<-\alpha<0$.  Moreover, each entry of $B^{V_0}(z(k))$ has type at most 
$2b_0$ and is bounded by $Ce^{2b_0|\Im k|}$.  Then
\begin{align*}\det(I+B^V(z(k)))& =\det (I+B^{V_0}(z(k))+D(z(k)))
\\ &  = \det(I+B^{V_0}(z(k)))
+ O\left(\frac{e^{|\Im k|(4b_0\#\{l:l\in \mce\}-\epsilon)}}{|k|}\right).
\end{align*}
Applying Theorem \ref{thm:complexanal} as in the proof of Theorems
\ref{thm:upperbound} and \ref{thm:asymptotics} finishes the proof.

It remains to show (\ref{eq:dbounds}).  Note that we may write, when
$|r_{j_0}(z)|$ is sufficiently large, 
\begin{multline*}
[I+(V_0+W)R_0(\omega_{\mce}(z))]^{-1}= \\
\left(I+\sum_{m=1}^{\infty}(-1)^m\left[ 
\left(I+V_0R_0(\omega_{\mce}(z))\right)^{-1}WR_0
(\omega_{\mce}(z))\right]^m\right)\left[I+V_0R_0(\omega_{\mce}(z))\right]^{-1}.
\end{multline*}
Then 
\begin{align*}
\varphi_{l+}^V(z) & = 
(I +VR_0(\omega_{\mce}(z)))^{-1}\left(V_0+W)\Phi_{l+}(\bullet,
z)\right)\\
 & = \left(I+V_0R_0(\omega_{\mce}(z))\right)^{-1}
\left(V _0 \Phi_{l+}( \bullet,
z)\right)
+ \left(I+V_0R_0(\omega_{\mce}(z))\right)^{-1}\left(W \Phi_{l\pm} ( \bullet,
z) \right)  \\
& \vspace{2mm}
 + \sum_{m=1}^{\infty} (-1)^m\left[ \left(I+V_0R_0(\omega_{\mce}(z))
\right)^{-1}WR_0
(\omega_{\mce}(z))\right]^m\left[I+V_0R_0(\omega_{\mce}(z))\right]^{-1}
\Phi_{l+}(\bullet,z).
\end{align*}
The first term on the right is $\varphi^{V_0}_{l+}(z)$.  The second term is,
as in (\ref{eq:mixedsignbd}),
bounded by 
$Ce^{(b_0-\epsilon)|\Im r_l(z)|}.$  Again as in (\ref{eq:mixedsignbd}),
 the third 
term is also bounded by $Ce^{(b_0-\epsilon)|\Im r_l(z)|}.$   Putting all this
into the definition of $b_{++lj}^{V}(z)$, we see that
$$b_{++lj}^{V}(z)= b_{++lj}^{V_0}(z) +
O\left(\frac{e^{(2b_0-\epsilon)|\Im r_{j_0}(z)|}}
{|r_{j_0}(z)|}\right).$$  A similar 
argument works for the other entries of $B^V(z)$, proving (\ref{eq:dbounds}).
\end{proof}

Combining the previous theorem with the results for potential scattering
in one dimension \cite{froese,zw1}, we obtain the following corollary.
\begin{cor}
Let $V(x,y)=V_0(x)+W(x,y) \in L^{\infty}_{\comp}(X;\Complex)$, where 
the support of $V_0$ is contained in $[-b,b]$ and in no smaller interval
and $\supp W\subset[-b+\epsilon,b-\epsilon]\times Y$ for some $\epsilon>0$.
Then on any sheet of $\hat{Z}$ with corresponding labeling set $\mce$,
\begin{multline*}
\#\{ z_k: z_k \; \text{is a pole of the resolvent of 
}\; \Delta +V \; \text{ on this sheet},\\
|r_{j_0}(z_k)|<r,\; \Im r_{j_0}(z_k)<-\alpha\} 
= \frac{4b}{\pi}\#\{l: \; l\in \tilde{\mce}\}  r+o_{\alpha}(r)
\end{multline*}
for any $\alpha>0$.
\end{cor}

\section{Results for half-cylinders} \label{s:halfcylinder}

In this section, we consider half-cylinders $X=[0,\infty)\times Y$, with 
$\Delta$ either the Dirichlet or Neumann Laplacian on $X$.  Let $V\in 
L^{\infty}_{\comp}(X;\Complex)$.  The resolvent $(\Delta +V -z)^{-1}$ has
a meromorphic continuation to $\hat{Z}$ just as in the full cylinder case.
We give several results analogous to the results for full cylinders.
Since the proofs are so similar, we only sketch them.

Let $R_{0\pm}(z)=(\Delta -z)^{-1}$ be the resolvent for the Neumann ($+$)
or Dirichlet ($-$) Laplacian on $X$, for $z\in \hat{Z}$.  Restrict $z$ to 
one fixed sheet of $\hat{Z}$, with corresponding labeling set 
$\mce $.  Then, following
the same argument as in the beginning of Section \ref{s:pfthmub},
we can show that when $|\Im \proj(\omega_{\mce}(z))|>\| \Im V\|_{L^{\infty}}$,
the poles of the resolvent of $\Delta +V$ correspond to the zeros of 
$\det(I+B_{\pm}(z))$.  Here we are again using
``$+$'' for the Neumann Laplacian and ``$-$'' for the Dirichlet Laplacian.
  To define $B_{\pm}(z)$, let 
\begin{align*}
\Phi_{\pm l}(x,y,z)& =(e^{i\tilde{r}_l(z)x}\pm e^{-i\tilde{r}_l(z)x})
\phi_l(y)\\
\varphi_{\pm l}(x,y,z)& = \left( (I+VR_{0\pm}(\omega_{\mce}(z)))^{-1}
(V\Phi_{\pm l})(\bullet,z)\right)
(x,y).
\end{align*}
Then $B_{\pm}(z)=(b_{\pm jk}(z))_{jk\in \tilde{\mce}}$, with
$$b_{\pm jk}(z)= \frac{i}{2\tilde{r}_j (z)} \int_X
\varphi_{\pm j}(x,y,z)\overline{\Phi}_{\pm k}(x,y,z) d\vol_X.$$

We obtain the following analog of Theorem \ref{thm:upperbound}.
\begin{thm} Let $X= [0,\infty)\times Y$ and 
let $V\in L^{\infty}_{\comp}(X;\Complex)$, with 
$\supp V \subset [0,b]\times Y$.
Fix a sheet of $\hat{Z}$, and suppose that $\Im r_{j_0}(z)<0$ on this
sheet.  
Then, there is a constant $c_{V, \mce}\geq 0$ such that
for any $\alpha>0$, 
\begin{multline*}
\#\{ z_k: z_k \text{is a pole of the resolvent on this sheet},
|r_{j_0}(z_k)|<r,\; \Im r_{j_0}(z_k)<-\alpha\} \\
= c_{V,\mce}r+o_{\alpha}(r)
\end{multline*}
The constant $c_{V,\mce}$ depends on the potential 
$V$ and the sheet.  Moreover, 
$$c_{V,\mce}
\leq \frac{2b}{\pi} 
\#\{l:\Im \tilde{r}_l(z)<0 \; \text {when }\; z \; 
\text{lies on this sheet}\}.$$
\end{thm}
\begin{proof}
Just as in the proof of Lemmas \ref{l:mixedbd} and
 \ref{l:samesignbd}, we can show that on our fixed sheet
$$|b_{\pm jk}(z)|\leq \frac{C}{|\tilde{r}_j(z)|}e^{2b|\Im \tilde{r}_k(z)|}.$$
Then the proof follows just as the proof of Theorem \ref{thm:upperbound}.
\end{proof}

\begin{thm} 
 Let $X=[0,\infty)\times Y$ and suppose that the support
of $V$ is contained in $[0,b]\times Y$ and the number $b$ cannot
be replaced by a smaller one.
Restrict ourselves
to a sheet of $\hat{Z}$ with $\Im r_j(z)<0$ if and only if $j=j_0$.
Suppose that $\nu_{j_0}^2$ is a simple eigenvalue of 
$\Delta_Y$, with $\nu_{j_0}^2=\sigma_{l_0}^2.$
 Suppose, in addition,
that 
\begin{equation*}
C|V_{l_0l_0}(x)|=C|\int_{Y}V(x,y)|\phi_{l_0}(y)|^2 d\vol_Y|
\geq |V(x,y)|,\; \text{for}\; |x-b|<\epsilon
\end{equation*}
for some $C,\epsilon>0$.
Then,
for any $\alpha>0$, 
\begin{multline*}
\#\{ z_k: z_k\; \text{is a pole of the resolvent on this sheet},
|r_{j_0}(z_k)|<r,\; \Im r_{j_0}(z_k)<-\alpha\} \\
= \frac{2}{\pi}b 
r +o_{\alpha}(r).
\end{multline*}
\end{thm}
\begin{proof}
In this case $B_{\pm}(z)$ is a single function, $b_{\pm l_0l_0}$.  Let
$k=\tilde{r}_{l_0}(z)$ and let $z(k)$ be the corresponding point on 
$\hat{Z}$.  
We have
\begin{align*}
b_{\pm l_0l_0}(z(k))& = \frac{i}{2k}\int_X 
(e^{ikx}\pm e^{-ikx})\overline{\phi}_{l_0}
\left[[I+VR_0(w_{\mce}(z(k)))]^{-1}V \Phi_{l\pm}
(\bullet,z(k))\right]d\vol_X\\
& = \frac{i}{2k}\int_X 
e^{ikx}\overline{\phi}_{l_0}\left[I+VR_0(w_{\mce}(z(k)))\right]^{-1}V 
f_{l_0}
(\bullet,z(k))  d\vol_X + O(e^{b|\Im k|}).
\end{align*}
Here $f_{l_0}(x,y,z)=e^{i\tilde{r}_{l_0}(z)x}\phi_{l_0}(y)$, and we have
used a bound similar to that of Lemma \ref{l:r0bound} to obtain the 
bound $O(e^{b|\Im k|})$ on the rest.  Following the technique of Lemmas
\ref{l:mixedbd} and \ref{l:mintype} shows that $b_{\pm l_0l_0}(z(k))$ is
an exponential function of type $2b$ for $\Im k \leq -\alpha$.  The proof is
completed as in the proof of Theorem
\ref{thm:asymptotics}.
\end{proof}

\small
\noindent
{\sc 
Department of Mathematics,
University of Missouri,
Columbia, Missouri 65211\\}
\end{document}